\documentclass[10pt,twocolumn]{article}

\usepackage{algorithm, algorithmic, setspace}
\usepackage{graphicx} % for pdf, bitmapped graphics files
\usepackage{amsmath} % assumes amsmath package installed
\usepackage{amssymb}  % assumes amsmath package installed
\usepackage{amsthm}
\usepackage{amsfonts}
\usepackage{color}
\usepackage{subfig}
\usepackage[english]{babel}
\usepackage[latin1]{inputenc}

\newcommand{\soft}{\mathbb{S}}

\newcommand{\X}{\mathcal{X}}

\newcommand{\argmin}[1]{\underset{#1}{\mathrm{argmin\,}}}

\newcommand{\reg}{\mathrm{\mathbf{Reg}}^d_T}

\newcommand{\R}{\mathbb{R}}
\newcommand{\N}{\mathbb{N}}

\newtheorem{lemma}{Lemma}

\newtheorem{theorem}{Theorem}
\newtheorem{corollary}{Corollary}
\newtheorem{remark}{Remark}

\newtheorem{assumption}{Assumption}

\title{\LARGE \bf Online optimization in dynamic environments:\\ a regret analysis for sparse problems}

\author{Sophie M. Fosson\\{\small{sophie.fosson@polito.it}}
}

\begin{document}

\maketitle
\thispagestyle{empty}
\pagestyle{empty}

%%%%%%%%%%%%%%%%%%%%%%%%%%%%%%%%%%%%%%%%%%%%%%%%%%%%%%%%%%%%%%%%%%%%%%%%%%%%%%%%
\begin{abstract}
Time-varying systems are a challenge in many scientific and engineering areas. Usually, estimation of time-varying parameters or signals  must be performed online, which calls for the development of responsive online algorithms. In this paper, we consider this problem in the context of the sparse optimization; specifically, we consider the  Elastic-net model, which promotes parsimonious solutions. Following the rationale in \cite{mok16}, we propose an online algorithm and we theoretically prove that it is successful in terms of dynamic regret. We then show an application to the problem of recursive identification of time-varying autoregressive models, in the case when the number of parameters to be estimated is unknown. Numerical results show the practical efficiency of the proposed method.
\end{abstract}

\section{INTRODUCTION}
In machine learning, signal processing, and control, many problems are innately time-varying. In dynamical systems, the challenge is to deal with  time-varying parameters, which has motivated a long standing research on online identification \cite{pra16, bra06, bra16, li11, bit94}. In signal processing, instead, the problem is to track time-varying signals (for example, moving targets). More in general, we talk about online optimization when the problem can be formulated as the minimization of a time-varying cost functional. In the last years, the literature on this topic has grown rapidly, with particular attention to the convex case. \emph{Online convex optimization} (OCO, \cite{sha12book,haz16book}) can be described as a game in which at each time step $t=1,\dots,T$ a player has to minimize a convex cost functional $f_{t}$ revealed by an adversary. In principle, the player might find  the minimum by convex programming, but the need for a real time response typically prevents it, which leads to the development of suboptimal strategies able to track the desired time-varying target. In some cases, a known dynamic is envisaged in the model, while in others a completely adversarial approach is considered, with no prior information on the evolution.

In the OCO literature, a benchmark to evaluate a strategy is provided by the regret \cite{sha12book,haz16book,zin03,haz07,mok16,hos16,sha18}, which basically measures the difference between the player's sequence of decisions and the best strategy in hindsight (\emph{i.e.,} the minimization of each $f_t$). A strategy is defined successful if its regret is sublinear in $T$. In the last decade, much attention has been devoted to the regret analysis in OCO: in \cite{zin03}, a gradient-based algorithm was proposed that achieves $\mathcal{O}(\sqrt{T})$ behavior; some years later, \cite{haz07} obtained $\mathcal{O}(\log{T})$ in the case of strongly convex functionals. This result was improved in \cite{mok16}, which obtained  $\mathcal{O}(1+C_T)$, where $C_T$ is the sum of the distances between successive reference points. Very recently, regret has been investigated also in OCO distributed settings \cite{hos16,sha18}.

In most of the mentioned works, $f_t$ is assumed to be differentiable, with bounded gradient. This does not envisage sparsity-promoting convex cost functionals, such as, Lasso \cite{tib96} and Elastic-net \cite{zou05}, which contain a (non-differentiable) $\ell_1$-norm. Such functionals are often used for model selection and data compression in estimation/identification problems \cite{hof13,pil14,yil18}. The key idea is that the $\ell_1$ regularization promotes sparse solutions (namely, solutions with many zeros), and can be used to build parsimonious models. This is of particular interest when models are built from large and noisy datasets. Moreover, in signal processing, sparsity-promoting functionals have gained new popularity with the advent of Compressed Sensing (CS, \cite{don06,fou13}), which states that sparse signals can be recovered from compressed linear measurements.

The aim of this work is to tackle the problem of sparse OCO. Since the best results in terms of regret analysis have been obtained so far with strongly convex functionals, we consider the (strongly-convex) Elastic-net model. We then elaborate an online strategy, based on iterative soft thresholding (IST, \cite{dau04, for10}), and we analyze its regret. Finally, we present a practical application for the identification of time-varying systems.

In the context of sparse signal processing and CS, the problem of tracking time-varying sparse signals has been tackled, but not in the regret analysis perspective. Moreover, in most works a known dynamics was assumed. In \cite{zin13}, an approximate message passing method was proposed to track time-varying sparse signals (acquired according to the CS paradigm), knowing that the dynamics is ruled by a Markov model; simulations showed good efficiency of the approximate message passing  in different applications, but no theoretical analysis was provided. In \cite{cha16}, a dynamic filtering via $\ell_{1}$ minimization was proposed, again assuming a Markov evolution model. The related optimization problem was formulated, at time $t$, as a Lasso plus a term promoting the consistency with the prediction of previous time step. Two algorithms were proposed: BPDN-DF, based on basis pursuit and Kalman filtering optimization, has theoretical guarantees on the error boundedness \cite[Theorem III.1]{cha16}, but degrades  when the  inaccuracy  in  the dynamics model increases; RWL1-DF is a re-weighted $\ell_1$ filtering, which achieves better performance according to numerical tests. \cite{mot15} considered the same model as \cite{cha16} (except for an $\ell_1$-norm term that enforces consistency with the prediction), then proposed an iterative algorithm and conditions for exact recovery.
A Kalman filtering approach was considered also in \cite{zac12}, where a dynamic iterative pursuit was proposed, based on a variant of orthogonal matching pursuit including prediction. In \cite{bal15}, IST was analyzed for tracking time-varying sparse signals, and a theoretical error bound was provided under the assumption of boundedness of the sparse signal and of its derivative.
%
%Tolgo questo: da discutere poi in una versione estesa
We finally mention the slightly different problem of streaming sets measurements, in which the unknown signal is static, but observations come sequentially, and one aims to refine the estimate online. Homotopy methods \cite{sal10,hof13} were proposed to tackle it. The approach is recursive: recovery is not performed from scratch each time, but past estimation is used to update the current one.

The paper is organized as follows. In Section \ref{sec:ps}, we state the optimization problem in a rigorous way. In Section \ref{sec:pm}, we introduce the proposed algorithm, whose regret is analyzed in Section \ref{sec:ra}. We then show a few numerical results in an online identification example (Section \ref{sec:nr}), and we finally draw some conclusions.
\section{PROBLEM STATEMENT}\label{sec:ps}
In OCO, at each time $t$ a convex cost functional $f_t$ is revealed to a player, that plays its action $x_{t+1}$. A suitable performance metric to evaluate the strategy is the so-called dynamic
regret, which is defined as follows \cite{mok16}: 
\begin{equation*}
\reg (z_1,\dots,z_T):=\sum_{t=1}^{T}{f_{t}(x_{t})-f_{t}(z_t)}
\end{equation*}
where 
\begin{equation}\label{def:zt}
z_t:=\argmin{x\in\X}f_{t}(x).
\end{equation}
where $\X\subseteq \R^n$ is a suitable state space.
We specify that the notion of static regret is considered in some works, where the desired target is the minimum of $\sum_t f_t$ \cite{mok16}. Dynamic regret instead is focused on tracking, which is the case of our interest.
We now formulate the problem using a sparse signal processing notation (the analogy with system identification will be clear in Section \ref{sec:nr}). We consider a time-varying sparse signal $v_t\in\R^n$, $t=1,\dots, T$, $T\in\N$, \emph{i.e.}, at each time step $t$, $v_t$ has $k_t\ll n$ non-zero components. We assume that linear measurements 
\begin{equation}\label{acqisition}
y_t=A_t v_t+e_t,~~A_t\in\R^{m, n}
 \end{equation}
are available, where $e_t$ is a possible measurement noise. Underdetermined problems with  $m<n$ (as in the CS setting) are envisaged. 
It is well known that the problem of sparse signal recovery can be recast into a convex optimization problem, leveraging the fact that $\ell_1$ norm promotes sparsity \cite{for10}. Specifically, minimizing $\frac{1}{2}\left\|y_t-A_t x\right\|_2^2+\lambda\left\|x\right\|_1$, $\lambda>0$ produces sparse solutions that are consistent with measurements: this is the popular Lasso problem \cite{tib96}. 

In this paper, following the rationale of \cite{mok16}, we propose to use the (strongly convex) Elastic-net model \cite{uma10}, which adds a Tikhonov regularization term to the Lasso and reads as follows:
\begin{equation}\label{eq:elasticnet}
	\begin{split}
	&f_t(x):=\frac{1}{2}\left\|y_t-A_t x\right\|_2^2+\lambda\left\|x\right\|_1+\frac{\mu}{2}\left\|x\right\|_2^2\\
	&\lambda>0, \mu>0, t=1,\dots,T.\\
	\end{split}
\end{equation}
A preliminary version of this model (with constant $A_t=A$) was introduced in \cite{fox17}. From a practical viewpoint, Tikhonov regularization promotes a grouping effect of correlated variables and is preferred in those sparse applications where it is undesired to discard correlated variables. In other terms, if two columns of $A_t$ are strongly correlated, Lasso would discard one of them, while Elastic-net would preserve both. For simplicity, we assume that the weight $\lambda$ is constant in $t$, which intrinsically assumes that the sparsity level $k_t$ is not significantly changing in time. 

Our goal is to solve the following optimization problem:
\begin{equation}\label{eq:theproblem}
\text{For any } t=1,\dots,T,~~~ \min_{x\in\R^n} f_t(x). 
\end{equation}
Since each $f_t$ is strongly convex, in principle a convex optimization algorithm  can be performed at each time step to get the desired minimum. Nevertheless, in the practice such algorithms are too slow to provide a solution in real time, in particular in case of large  $n$. For this motivation, we propose to investigate a low-complexity iterative method that can online estimate $v_t$. Such method is expected to be suboptimal, but will turn out to be successful in terms of dynamic regret and efficient in numerical experiments. 
\section{PROPOSED METHOD}\label{sec:pm}
Our idea is to adapt IST for problem \eqref{eq:theproblem} and analyze its dynamic regret. IST is a gradient-based algorithm which has been proved to converge to a minimum of Lasso \cite{for10} and to the minimum of the Elastic-net \cite{uma10}. Our goal is to prove that it can be used successfully also in the dynamic case. 

The choice of IST has different motivations. On the one hand, its simplicity allows a straightforward implementation and makes its theoretical analysis affordable. On the other hand, IST is prone to decentralization and parallel implementation \cite{rfm15,fox14,fox16,fia13}. In this work, we deal with a centralized setting, but the extension to a distributed setting (for example, a sensor network) might be investigated in the future. A drawback of IST is its convergence slowness, in terms of number of iterations, which makes other iterative algorithms more popular, firstly the alternating direction method of multipliers (ADMM, \cite{boy10}). The implementation and the analysis proposed in this paper for IST can be extended to ADMM, which will be object of future work.

IST for (static) Elastic-net was derived in \cite{uma10}. Since our derivation is slightly different, we now describe it in detail.
\subsection{Batch IST for Elastic-net}
Let us define:
\begin{equation*}
\begin{split}
&g(x):=\frac{1}{2}\left\| y-Ax\right\| _{2}^{2}+\frac{\mu}{2}\|x\|_{2}^{2};~~f(x):=g(x)+\lambda\|x\|_{1}\\
&s(x,b):=\frac{1}{2\tau}\left\| x-b\right\| _{2}^{2}-\frac{1}{2}\left\| Ax-Ab\right\| _{2}^{2}\\
&f(x,b):=f(x)+s(x,b)\\
\end{split}
\end{equation*}
where $x,b\in\R^{n}$, $A\in\R^{m, n}$, $y\in\R^m$, and $\tau>0$ is chosen so that $\tau\|A\|_{2}^{2}<1$; in this way, $s(x,b)\geq0$.
As a consequence, following \cite{for10}, we can easily see that minimizing $f(x)$ is equivalent to minimizing the so-called surrogate
functional $f(x,b).$ In particular, if $z$ is the minimum of $f(x)$, then $(z,z)$ is a minimum for $f(x,b)$. 

Adding $s(x,b)$, we cancel the term $\|Ax\|_{2}^{2}$, then, when $b$ is fixed the problem $\min_{x\in\R^{n}}f(x,b)$ is separable in the single components of $x$ and its solution can be written in closed form. In fact, the terms of $f(x,b)$ depending on $x$ can be grouped as follows:
\begin{align*}
& \frac{1}{2}\left(\mu+\frac{1}{\tau}\right)\|x\|_{2}^{2}-\langle x,\frac{b}{\tau}+A^{T}(y-Ab)\rangle+\lambda\|x\|_{1}=\\
& =\frac{1}{2}\left(\mu+\frac{1}{\tau}\right)\left\|x-\frac{b+\tau A^{T}(y-Ab)}{1+\mu\tau}\right\|_{2}^{2}+\lambda\|x\|_{1}+c
\end{align*}
where $c\in\R$ does not depend on $x$.
We then conclude:
\begin{equation}
\argmin{x\in\R^{n}}f(x,b)=\soft_{\frac{\lambda\tau}{\mu\tau+1}}\left[\frac{b+\tau A^{T}(y-Ab)}{1+\mu\tau}\right]
\end{equation}
where $\soft_{\beta}:\R^{n}\to\R^{n}$, $\beta>0$, is the well-known component-wise soft thresholding operator \cite{for10}, defined as follows: for $z\in\R$, $\soft_{\beta}[z]=z-\beta$ if $z>\beta$; $\soft_{\beta}[z]=z+\beta$ if $z<-\beta$; $\soft_{\beta}[z]=0$ otherwise.

On the other hand, fixed $x$, it is easy to check that $x=\argmin{b\in\R^{n}}f(x,b)$.
Alternating the minimizations in $x$ and $b$, we obtain Algorithm \ref{alg:IST_static}, which converges to the minimum of $f(x)$ (this can be easily proved exploiting the contractivity). Practically, the algorithm stops when a suitable numerical convergence criterion is met. Notice that in this batch procedure the time $t$ only refers to the algorithm's iterations.
\begin{algorithm}
	\caption{Batch IST for Elastic-net}
	\label{alg:IST_static}
	\begin{algorithmic}[1] 
		\STATE $x_0=0$
		\FOR{$t=1,\dots,T_{stop}$}
		\STATE $x_t=\soft_{\frac{\lambda\tau}{\mu\tau+1}}\left[\frac{x_{t-1}+\tau A^{T}(y-Ax_{t-1})}{1+\mu\tau}\right]$
		\ENDFOR
	\end{algorithmic}
\end{algorithm} 
\subsection{Online IST for dynamic Elastic-net}
In order to tackle the dynamic problem \eqref{eq:elasticnet}-\eqref{eq:theproblem}, we propose an online version of Algorithm \ref{alg:IST_static}, which performs an IST step at each time $t$. This is summarized in Algorithm \ref{alg:IST_dynamic}, which will be theoretically analyzed in Section \ref{sec:ra}. If the computational resources allow it, one could perform $r>1$ IST steps at each $t$, as we will discuss in Section \ref{sec:nr}.
\begin{algorithm}
	\caption{Online IST for dynamic Elastic-net}
	\label{alg:IST_dynamic}
	\begin{algorithmic}[1] 
		\STATE $x_0=0$
		\FOR{$t=1,\dots,T$}
		\STATE $x_t=\soft_{\frac{\lambda\tau}{\mu\tau+1}}\left[\frac{x_{t-1}+\tau A_t^{T}(y_t-A_tx_{t-1})}{1+\mu\tau}\right]$
		\ENDFOR
	\end{algorithmic}
\end{algorithm}
\subsection{Related literature}
We observe that online IST is an instance of Composite Objective Mirror Descent (COMID, \cite{duc10}), a popular technique in machine learning. Developed in the same years in different frameworks, both IST and COMID come from the forward-backward techniques. COMID was introduced to tackle online optimization with a regularizer term, \emph{i.e.}, the minimization of cost functionals of kind $\sum_{t}f_{t}(x)+r(x)$, where $r(x)$ is a (static) regularizer. $f_{t}$ and $r$ are both supposed to be convex. Its update is given by $x_{t+1}=\argmin x\langle f'_{t}(x_{t}),x\rangle+B_{\psi}(x,x_{t})+\eta r(x)$, 
where $B_{\psi}$ is a Bregman divergence. In \cite{duc10}, the static regret was proved to behave as $\mathcal{O}(\sqrt{T})$ and $\mathcal{O}(\log T)$,
for convex and strongly convex functions, respectively. The assumptions
to obtain this result, however, include the boundedness of $\|f'_{t}\|$ and
$\eta$ vanishing as $\frac{1}{\sqrt{T}}.$ These assumptions are not required in this paper. In particular, vanishing parameters are generally not desirable when the time horizon $T$ is large or tends to infinity. 
\section{REGRET ANALYSIS}\label{sec:ra}
In this section, we prove that $\reg$ for Online IST for dynamic Elastic-net is sublinear; in particular, it only depends on the evolution of $z_t$ \eqref{def:zt} (or similarly, on the evolution of $A_t$ and $v_t$ \eqref{acqisition}). We now provide some lemmas that build the proof of our main result. At the end of the section, we discuss some consequences of it. In the following, we assume that the evolution of the system is bounded.
\begin{assumption}
For $t=1,\dots,T$, $\tau \|A_t\|_2^2\leq 1$, and $\|v_t\|_2\leq v_M$ for some $v_M>0$. 
\end{assumption}
As a consequence, also $\|z_t\|_2<z_M$ for some $z_M>0$. We remark that this an assumption on the system's evolution, while we do not force any boundedness on the algorithm's evolution $x_t$ or $f_t(x_t)$. This is an improvement with respect to \cite{mok16}, where the boundedness of $\nabla f_t$ was required \cite[Assumption 3]{mok16}, which for example is not satisfied for quadratic cost functions, unless a bounded state space is assumed \cite[Example in Section V]{mok16}.

Given any $x\in\R^n$, let us define 
\begin{equation*}
\Gamma_t(x):=\soft_{\frac{\lambda\tau}{\mu\tau+1}}\left[\frac{x+\tau A_t^{T}(y_t-A_tx)}{1+\mu\tau}\right].
\end{equation*}
\begin{lemma}\label{lem:contr}[Contractivity] Let $z_t=\argmin{x\in\R^n} f_t(x)$. Then, for any $x\in\R^n$, 
$$\|\Gamma_t(x)-z_t\|_2\leq \frac{1}{1+\mu\tau}\|x-z_t\|_2.$$
\end{lemma}
\begin{proof}
Since $z_t$ is the minimum, it is also a stationary point: $\Gamma_t(z_t)=z_t$. Since $\soft_{\beta}$ is non-expansive, we have:
\begin{align*}
%\|\Gamma_t(x)-&z_t\|_2\leq\frac{\left\|I-\tau A_t^T A_t \right\|_2\|x-z_t\|_2}{1+\mu\tau}=\frac{\|x-z_t\|_2}{1+\mu\tau}
&\|\Gamma_t(x)-z_t\|_2 =\|\Gamma_t(x)-\Gamma_t(z_t)\|_2\\
& \leq\left\|\frac{x+\tau A_t^{T}(y_t-A_t x)}{1+\mu\tau}-\frac{z_t+\tau A_t^{T}(y_t-A_tz)}{1+\mu\tau}\right\|_2\\
& \leq\frac{\left\|I-\tau A_t^T A_t \right\|_2}{1+\mu\tau}\|x-z_t\|_2=\frac{1}{1+\mu\tau}\|x-z_t\|_2
\end{align*}
where $I\in\R^{n, n}$ is the identity matrix.
\end{proof}
Using the contractivity, we can prove that, at a given $t$, the distance between the played action and the current minimum is controlled by the distance between successive minima. Let $$\Delta_t:=\|z_t-z_{t-1}\|_2.$$
\begin{lemma}\label{lem:somme} For any $t=1,\dots,T$,
	\begin{equation*}
	\begin{split}
	(a)~~&\sum_{t=2}^T \left\|x_t -z_t\right\|_2\leq c_1 + c_2 \sum_{t=2}^T \Delta_t\\
	(b)~~&\sum_{t=2}^T \left\|x_t -z_t\right\|^2\leq c_3 + c_4 \sum_{t=2}^T \Delta_t^2+c_5 \sum_{t=2}^T \Delta_t\\
	\end{split}
	\end{equation*}
	where $c_i>0$, $i=1,\dots,5$ are assessed in the proof.
\end{lemma}
\begin{proof}
Using the triangle inequality and Lemma \ref{lem:contr}, for any $t=2,\dots,T$,
%$\|x_t-z_t\|_2\leq \frac{1}{1+\mu\tau}\|x_{t-1}-z_{t-1}\|_2 +\Delta_t.$
\begin{equation*}
\begin{split}
\|x_t-z_t\|_2&\leq \|x_t-z_{t-1}\|_2 +\Delta_t\\
&\leq \frac{1}{1+\mu\tau}\|x_{t-1}-z_{t-1}\|_2 +\Delta_t.\\
\end{split}
\end{equation*}
Summing over $t=2,\dots,T$, 
%$\frac{\mu\tau}{1+\mu\tau}\sum_{t=2}^T \|x_t-z_t\|_2\leq  \frac{1}{1+\mu\tau}\left(\|x_{1}-z_{1}\|_2 - \|x_{T}-z_{T}\|_2\right) + \sum_{t=2}^T \Delta_t.$
 \begin{equation*}
 \begin{split}
 &\frac{\mu\tau}{1+\mu\tau}\sum_{t=2}^T \|x_t-z_t\|_2\leq  \frac{\|x_{1}-z_{1}\|_2 - \|x_{T}-z_{T}\|_2}{1+\mu\tau}+\\&~~~ + \sum_{t=2}^T \Delta_t.
 \end{split}
 \end{equation*}
We have then proved $(a)$ with 
%$c_1= \frac{1}{\mu\tau}\left(\|x_{1}-z_{1}\|_2 - \|x_{T}-z_{T}\|_2\right)$ and $c_2=\frac{1+\mu\tau}{\mu\tau}.$
\begin{equation*}
\begin{split}
c_1&= \frac{1}{\mu\tau}\left(\|x_{1}-z_{1}\|_2 - \|x_{T}-z_{T}\|_2\right);\\
c_2&=\frac{1+\mu\tau}{\mu\tau}.\\
\end{split}
\end{equation*}
To prove $(b)$, we use: 
\begin{equation*}
\begin{split}
&\|x_t-z_t\|_2^2\leq \|x_t-z_{t-1}\|_2^2 +\Delta_t^2+2\|x_t-z_{t-1}\|_2\Delta_t \\
&\leq \frac{\|x_{t-1}-z_{t-1}\|_2^2}{(1+\mu\tau)^2} +\Delta_t^2+\frac{4 z_M}{1+\mu\tau}\|x_{t-1}-z_{t-1}\|_2.\\
\end{split}
\end{equation*}
Summing over $t=2,\dots,T$, we obtain $(b)$ with
%$c_3= \frac{1}{\mu\tau(\mu\tau+2)}\left(\|x_{1}-z_{1}\|_2^2 - \|x_{T}-z_{T}\|_2^2\right)+\frac{4 z_M}{1+\mu\tau}\left(\|x_{1}-z_{1}\|_2 - \|x_{T}-z_{T}\|_2\right)$, 
%$c_4=\frac{(1+\mu\tau^2)}{\mu\tau(\mu\tau+2)}$, $c_5=\frac{4 z_M}{1+\mu\tau}.$
\begin{equation*}
\begin{split}
c_3&= \frac{\|x_{1}-z_{1}\|_2^2 - \|x_{T}-z_{T}\|_2^2}{\mu\tau(\mu\tau+2)}+\\&~~~+\frac{4 z_M}{1+\mu\tau}\left(\|x_{1}-z_{1}\|_2 - \|x_{T}-z_{T}\|_2\right);\\
c_4&=\frac{(1+\mu\tau^2)}{\mu\tau(\mu\tau+2)};~~~c_5=\frac{4 z_M}{1+\mu\tau}.
\end{split}
\end{equation*}
\end{proof}
\begin{lemma}\label{lem:surro}
$$f_{t-1}(x_t)-f_{t-1}(z_{t-1})\leq \frac{1}{\tau}\left\| x_{t-1}-z_{t-1}\right\|_2^2.$$
\end{lemma}
\begin{proof} Notice that $$f_{t-1}(x_t) \leq f_{t-1}(x_t,x_{t-1}) \leq f_{t-1}(z_{t-1},x_{t-1}),$$ using the fact that $x_t$ minimizes $f(\cdot, x_{t-1})$. Now, %$f_{t-1}(z_{t-1},x_{t-1})-f_{t-1}(z_{t-1})=s(x_{t-1},z_{t-1})\leq\frac{1}{\tau}\left\| x_{t-1}-z_{t-1}\right\|_2^2.$ 
 \begin{equation*}
 \begin{split}
 f_{t-1}(z_{t-1},x_{t-1})-f_{t-1}(z_{t-1})&=s(x_{t-1},z_{t-1})\\
 &\leq\frac{\left\| x_{t-1}-z_{t-1}\right\|_2^2}{\tau}.
 \end{split}
 \end{equation*}
\end{proof}
\begin{lemma}\label{lem:d_t}
	Let 
	%$d_t(x):=f_t(x)-f_{t-1}(x).$ 
	\begin{equation*}
	\begin{split}d_t(x)&:=f_t(x)-f_{t-1}(x)\\&=\frac{1}{2}\left\|A_t x-y_t\right\|_2^2-\frac{1}{2}\left\|A_{t-1} x-y_{t-1}\right\|_2^2.
	\end{split}
	\end{equation*} 
	Then, %$d_t(x_t)-d_t(z_t)\leq \gamma_1 \left\| x_t-z_t\right\|_2 +\gamma_2 \left\| x_t-z_t\right\|_2^2,$
 	\begin{equation*}
 	\begin{split}
 	d_t(x_t)-d_t(z_t)\leq \gamma_1 \left\| x_t-z_t\right\|_2 +\gamma_2 \left\| x_t-z_t\right\|_2^2
 	\end{split}
 	\end{equation*} 
	where $\gamma_1>0, \gamma_2>0$ are assessed in the proof.
\end{lemma}
\begin{proof}
Observing that  $\left\| x_t+z_t\right\|_2\leq \left\| x_t-z_t\right\|_2+2\left\|z_t\right\|_2$, straightforward computation leads to the following bound:
\begin{equation*}
\begin{split}
&d_t(x_t)-d_t(z_t)\leq   \left\| x_t-z_t\right\|_2  \left\| A_t^T y_t -A_{t-1}^T y_{t-1}\right\|_2\\&~~+ \frac{1}{2}  \left\| x_t-z_t\right\|_2  \left\| A_t^T A_t-A_{t-1}^T A_{t-1}\right\|_2  \left\| x_t+z_t\right\|_2\\
& \leq   \frac{2 v_M + 2 z_M}{\tau}\left\| x_t-z_t\right\|_2 + \frac{1}{\tau} \left\| x_t-z_t\right\|_2^2.
\end{split}
\end{equation*}
%$d_t(x_t)-d_t(z_t)\leq \frac{2 v_M + 2 z_M}{\tau}\left\| x_t-z_t\right\|_2 + \frac{1}{\tau} \left\| x_t-z_t\right\|_2^2.$
\end{proof}
Based on these lemmas, we can prove our main result.
\begin{theorem}\label{theo:regret}
$$	\reg \leq \alpha_0+\alpha_1 \sum_{t=2}^T \left\|z_t-z_{t-1}\right\|_2 + \alpha_2 \sum_{t=2}^T \left\|z_t-z_{t-1}\right\|_2^2$$
where $\alpha_i>0$, $i=0,1,2$ are assessed in the proof.
\end{theorem}
\begin{proof}
First, we compute a bound for the loss:
\begin{equation*}
\begin{split}
&f_t(x_t)-f_t(z_t)=\\&=f_t(x_t)-f_t(z_t)\pm f_{t-1}(x_t)+ f_{t-1}(z_{t})- f_{t-1}(z_{t-1})\\
&\leq d_t(x_t)-d_t(z_t)+f_{t-1}(x_t)-f_{t-1}(z_{t-1})\\
\end{split}
\end{equation*}
%$f_t(x_t)-f_t(z_t)\leq d_t(x_t)-d_t(z_t)+f_{t-1}(x_t)-f_{t-1}(z_{t-1})$, 
where we have used the fact that $ f_{t-1}(z_{t})\leq f_{t-1}(z_{t-1})$. Now, applying Lemma \ref{lem:d_t} and Lemma \ref{lem:surro} we get: $f_t(x_t)-f_t(z_t)\leq  \gamma_1 \left\| x_t-z_t\right\| +\gamma_2 \left\| x_t-z_t\right\|_2^2+ \frac{1}{\tau}\left\| x_{t-1}-z_{t-1}\right\|_2^2.$
%\begin{equation*}
%	\begin{split}
%f_t(x_t)-f_t(z_t) &\leq  \gamma_1 \left\| x_t-z_t\right\| +\gamma_2 \left\| x_t-z_t\right\|_2^2+\\&+ \frac{\left\| x_{t-1}-%z_{t-1}\right\|_2^2}{\tau}.
%\end{split}
%\end{equation*}
Finally, summing over $t=2,\dots,T$ and exploiting Lemma \ref{lem:somme}, we obtain the thesis with
$$\alpha_0= \left(\frac{1}{\tau}+\gamma_2\right)+\gamma_1 c_1+c_3+\frac{1}{\tau}\left\|x_1-z_1\right\|_2^2-\frac{1}{\tau}\left\|x_T-z_T\right\|_2^2,$$ $$\alpha_1=\left(\frac{1}{\tau}+\gamma_2\right) c_5+\gamma_1 c_2,~~~\alpha_2=\left(\frac{1}{\tau}+\gamma_2\right)c_4.$$
\end{proof}
We notice that Theorem \ref{theo:regret} can be reformulated by highlighting that the dynamic regret is controlled by the evolution of the system variables $v_t, A_t$ (instead of $z_t$).
\begin{corollary}\label{cor:regret}
	$$\reg \leq  \alpha_0+\alpha_1 \sum_{t=2}^T \left\|\theta_t\right\| + \alpha_2 \sum_{t=2}^T \left\|\theta_t\right\|_2^2,$$
	where $\theta_t= \frac{(z_M+v_M) \left\| A_t^T A_t-A_{t-1}^T A_{t-1}\right\|_2}{\mu} + \frac{1}{\mu\tau}\left\|v_t-v_{t-1}\right\|_2.$
\end{corollary}
\begin{proof}
We prove that $\Delta_t=\|z_t-z_{t-1}\|_2\leq \theta_t$. Given this, the thesis follows from Theorem \ref{theo:regret}.
At any $t$, using the stationarity of $z_t$ and $z_{t-1}$ for $\Gamma_t$ and $\Gamma_{t-1}$, respectively: %$\Delta_t \leq \frac{\left\|(I-\tau A_t^t A_t)z_t-(I-\tau A_{t-1}^t A_{t-1})z_{t-1}\right\|_2}{1+\mu\tau}+\tau\frac{\left\|A_t^T A_t v_t-A_{t-1}^T A_{t-1} v_{t-1} \right\|_2}{1+\mu\tau}.$
 \begin{equation*}
 \begin{split}
 &\Delta_t \leq \frac{\left\|(I-\tau A_t^t A_t)z_t-(I-\tau A_{t-1}^t A_{t-1})z_{t-1}\right\|_2}{1+\mu\tau}\\&~~+\tau\frac{\left\|A_t^T A_t v_t-A_{t-1}^T A_{t-1} u_{t-1} \right\|_2}{1+\mu\tau}.\\
 \end{split}
 \end{equation*}
Adding and subtracting $(I-\tau A_{t-1}^t A_{t-1})z_{t}$ and $ \tau A_{t-1}^T A_{t-1}v_t$ within the two norms, respectively, we obtain:
%$\Delta_t \leq  \frac{(z_M+v_M) \tau \left\| A_t^T A_t-A_{t-1}^T A_{t-1}\right\|_2}{1+\mu\tau}+ \frac{\left\|z_t- z_{t-1}\right\|_2+\left\|v_t- u_{t-1}\right\|_2}{1+\mu\tau}.$
\begin{equation*}
\begin{split}
&\Delta_t \leq  \frac{(z_M+v_M) \tau \left\| A_t^T A_t-A_{t-1}^T A_{t-1}\right\|_2}{1+\mu\tau}\\&~~ + \frac{\left\|z_t- z_{t-1}\right\|_2+\left\|v_t- u_{t-1}\right\|_2}{1+\mu\tau}.\\
\end{split}
\end{equation*}
Therefore, %$\Delta_t \leq  \frac{(z_M+v_M) \left\| A_t^T A_t-A_{t-1}^T A_{t-1}\right\|_2}{\mu} + \frac{\left\|v_t- u_{t-1}\right\|_2}{\mu\tau}.$
 \begin{equation*}
 \begin{split}
 &\Delta_t \leq  \frac{(z_M+v_M) \left\| A_t^T A_t-A_{t-1}^T A_{t-1}\right\|_2}{\mu}\\& + \frac{\left\|v_t- v_{t-1}\right\|_2}{\mu\tau}.
 \end{split}
 \end{equation*}

\end{proof}
\begin{remark}
The regret analysis can be performed in an analogous way if $r>1$ IST steps are played at each time step. Starting from Lemma \ref{lem:contr}, we would have a stronger contractivity with contraction factor $(1+\mu\tau)^{-r}$, which would lead step by step to tighter constant parameters in Theorem \ref{theo:regret}. Since the purpose of this work is to describe the general behavior, a precise analysis of the improvement obtained with $r>1$ is left to future work.
\end{remark}
\subsection{Consequences}
Theorem \ref{theo:regret} and Corollary \ref{cor:regret} state a regret bound depending uniquely on the behavior of the system (and not directly on $T$) in terms of sum of (linear and quadratic) distances between successive iterates (in terms of minima $z_t$, or $A_t$ and $v_t$). This is in line with the results of \cite{mok16}, which obtained similar (linear) bounds for generic strongly convex problems, under stronger assumptions. In particular, as in \cite{mok16}, we can observe that: 
\begin{itemize}
 \item[-] if $z_t=z$ (or $A_t=A$ and $v_t=v$) for any $t>t_0$, then $\reg=\mathcal{O}(1)$, which implies that $x_t$ converges to  $z$; 
 \item[-] if $\|z_t-z_{t-1}\|_2=\mathcal{O}\left(\frac{1}{t^\eta}\right)$ (or $\|A_t^T A_t-A_{t-1}^T A_{t-1}\|_2=\mathcal{O}\left(\frac{1}{t^\eta}\right)$ and $\|v_t-v_{t-1}\|_2=\mathcal{O}\left(\frac{1}{t^\eta}\right)$), with $\eta>0$, then $x_t$ converges to $z_t$; 
 \item[-] if $\|z_t-z_{t-1}\|_2\leq C$ (or $\|A_t^T A_t-A_{t-1}^T A_{t-1}\|_2\leq C$ and $\|v_t-v_{t-1}\|_2\leq C$) where $C$ is constant, we have a steady state estimation error.
\end{itemize}
\section{NUMERICAL RESULTS}\label{sec:nr}
\begin{figure*}[ht]
	\centering
	%\begin{subfigure}[b]{0.25\textwidth}
	\includegraphics[width=0.98\columnwidth]{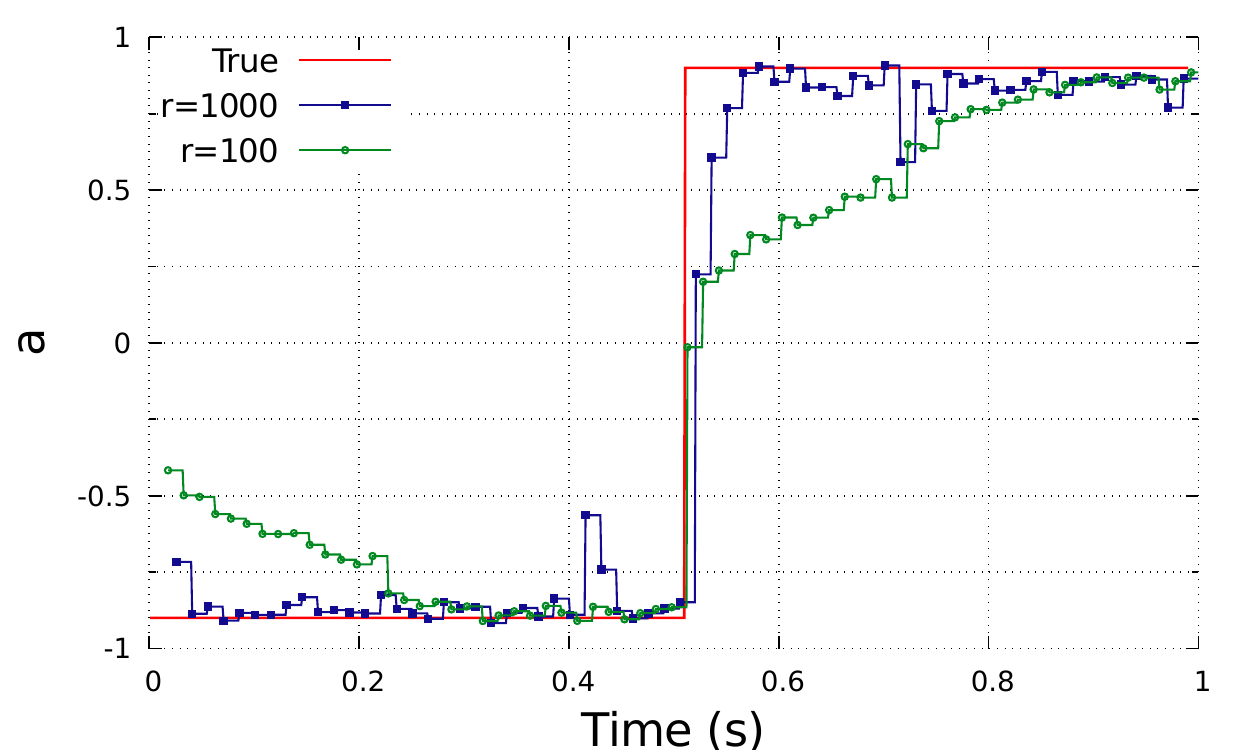}\quad
	\includegraphics[width=0.98\columnwidth]{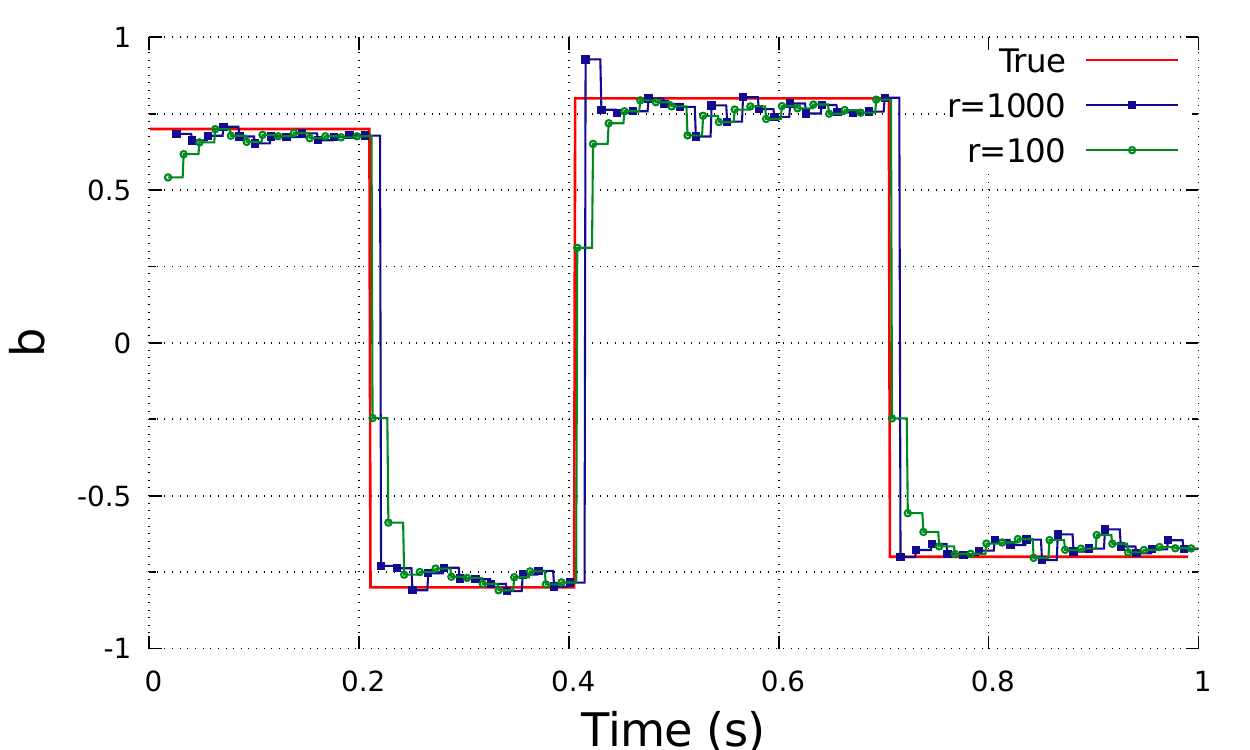}
	\caption{Estimation of $a_{1,t}$ and $b_{1,t}$, with online IST with $r=1000$ and $r=100$.}
	\label{fig:simulations}
\end{figure*}
\begin{figure}[ht]
	\centering
	%\begin{subfigure}[b]{0.25\textwidth}
	\includegraphics[width=0.98\columnwidth]{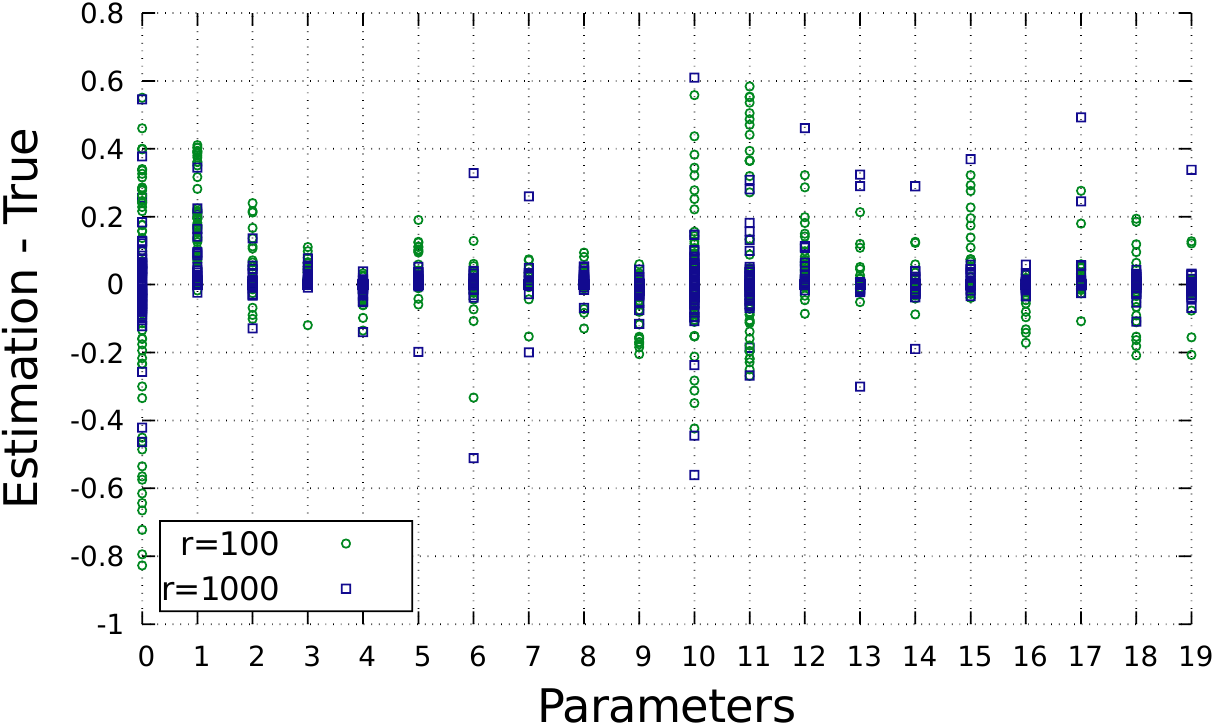}
	\caption{Signed error $\widehat{v}_{i,sm}-v_{i,sm}$ ($m=15$, $s=1,\dots,66$) for the $P+Q=20$ possible parameters (0 and 10 respectively are $a_{1,t}$ and $b_{1,t}$, while the others are the null parameters).}
	\label{fig:support}
\end{figure}
In this section, we present the implementation of the proposed method to tackle the online estimation of time-varying parameters in a time-varying autoregressive model with an exogenous input (TVARX, \cite{li11}). In a TVARX model, the input-output relationship is as follows: $y_t=\sum_{p=1}^P a_{p,t}y_{t-p}+\sum_{q=1}^Q b_{q,t}u_{t-q}+e_t$, 
where $u_t,y_t\in\R$ respectively are the measurable input and output; $e_t\in\R$ is the measurement error; $a_{p,t}, b_{q,t}\in\R$ are the time-varying parameters to be estimated. The goal is to perform a recursive identification, that is, to estimate $a_{p,t}$ and $b_{q,t}$ at time $t$, exploiting the knowledge of $u$ and $y$ and an estimation at time $t-1$. In our simulations, we assume that the dimensions $P$ and $Q$ are not known. A suitable strategy to overcome this lack is to fix a sufficiently large upper bound for them, and then look for a sparse solution, that is, to the most parsimonious model that represents the system accurately. 

The model and methodology presented in Sections \ref{sec:ps} and \ref{sec:pm} perfectly match this purpose. More precisely, we consider the following setting: we iteratively collect groups of $m$ output measurements $\mathbf{y}_t:=(y_t,\dots, y_{t+m})^T$, and we run $r$ steps of the online IST algorithm. It is easy to check that we can define the time-varying matrix $A_t\in\R^{m,P+Q}$ as follows:
\begin{equation*}
\left(\begin{array}{cccccc}
		y_{t-1}&\cdots&y_{t-P}&u_{t-1}&\cdots&u_{t-Q}\\
		y_{t}&\cdots&y_{t-P+1}&u_{t}&\cdots&u_{t-Q+1}\\
		\vdots&&&&&\vdots\\
		y_{t+m-1}&\cdots&y_{t+m-P}&u_{t+m-1}&\cdots&u_{t+m-Q}\\
	\end{array}\right).
\end{equation*}
Assuming $m$ sufficiently small so that $a$ and $b$ are constant between $t$ and $t+m$, we  formulate the problem as follows: at each $t=sm$, $s\in\N$, we collect $m$ measurements $\mathbf{y}_t$ and given $A_t$ we aim to recover $$v_t=(a_{1,t},\dots, a_{P,t}, b_{1,t},\dots,b_{Q,t})^T$$ from $$\mathbf{y}_t= A_t v_t +\mathbf{e}_t$$
where $\mathbf{e}_t:=(e_t,\dots, e_{t+m})^T$. $v_t$ is assumed to be sparse with respect to overestimated dimension $P+Q$.
For our simulations, we retrieve the TVARX(1,1) example considered in \cite[Section V]{li11}: $y_t=a_{1,t}y_{t-1}+b_{1,t}u_{t-1}+e_t$, with $P=Q=1$. As mentioned, we assume $P$ and $Q$ unknown and we overestimate them as $P=Q=10$. In terms of sparse signal recovery, we can say that we have to online estimate a time-varying sparse signal $v_t\in\R^n$ with $n=20$ and (unknown) sparsity $k=2$. Such sparse signal has constant support (but this information is assumed to be unknown and then not exploited in our procedure). We consider a time horizon of 1 second and sampling frequency of 1000 Hz. We assume that $a_{1,t}$ and $b_{1,t}$ are step-wise constant, with some abrupt changes (with respect to the example in \cite{li11}, we consider larger gaps). Specifically, we set:
{\small{
\begin{equation*}
a_1(t)=\left\{\begin{split}
&-0.9 \text { if } t<0.5\\
&0.9 \text{ otherwise.}
\end{split}\right.
\end{equation*}
\begin{equation*}
b_1(t)=\left\{\begin{split}
&0.7 \text { if } t<0.2\\
&-0.8 \text { if } 0.2 \leq t<0.4\\
&0.8 \text { if } 0.4 \leq t<0.7\\
&-0.7 \text{ otherwise.}\\
\end{split}\right.
\end{equation*}
}}
The input $u$ is assumed to be a standard Gaussian sequence, periodic with period $m$ (so that the right part of $A_t$ is constant for any $t=ms$, $s\in\N$). The measurement noise is white Gaussian with SNR around 20dB.

Finally, we assume $m=15$, therefore, as in CS, we have less measurements than parameters to estimate. In  the literature, this is also known as Compressive System Identification \cite{tot11,san11}. In \cite{tot11}, the case of constant parameters was considered, while in \cite{san11} an extension to piecewise-constant parameters was proposed. $A_t$ is not a typical CS sensing matrix,  but the randomness introduced by the chosen  $u$ and the circulant structure are promising features to use it in a CS framework \cite{for10}. In \cite{tot11,san11}, the properties of $A_t$ for CS where theoretically studied, but the obtained bounds are not tight \cite[Section IV.C]{san11}. This point is then still open and will be considered for future work, as well as a thorough  comparison to the approach proposed in \cite{san11}.

Coming back to our experiment, we run online IST with the following design parameters: $\lambda=2\times 10^{-2}$, $\mu=10^{-6}$, $\tau=3\times 10^{-2}$; the initial condition is assumed to be zero. At each $t=ms$, we run $r$ IST steps, with $r=1000$ and $r=100$. Simulations have been performed on a 2.67 GHz CPU, where $r=1000$ and $r=100$  respectively require 10 ms and 2 ms. Since we consider blocks of $m=15$ measurements (acquired in 15 ms), in both cases the algorithm does not exceed the time of acquisition. The overall delay for acquiring the measurements and running the simulations is then around 25 ms and 17 ms, respectively. We specify that performing a complete IST at each block (and then get the minimum of $f_t$) is not feasible in real time, since the convergence on this problem requires up to 35000 iterations and 500 ms for execution on the considered CPU.

We have performed 250 runs, whose average results are shown in Table \ref{tab:0}. The mean square error is defined as  MSE=$\frac{1}{P+Q}\sum_{s=1}^{T/m}\sum_{i=1}^{P+Q}\left\|v_{i,sm}-\widehat{v}_{i,sm}\right\|_2^2$, where $v_t=(a_{1,t},\dots, a_{P,t}, b_{1,t},\dots,b_{Q,t})^T$ and $\widehat{v}_t$ is our estimation. As expected, we can notice an average improvement when $r$ increases.
\begin{table}
	\begin{center}
		\caption{Average MSE and standard deviation over 250 runs, with $r$ IST steps at each $sm$.}
		%\centering
		\resizebox{0.5\columnwidth}{!} {
			\begin{tabular}{| r | c | c |}
				\hline
				IST steps  & Average MSE  & Std \\
				\hline
				$r=100$ & 0.011 & 0.002\\
				$r=1000$ & 0.006& 0.001\\
				\hline
			\end{tabular} 
			\label{tab:0}}
	\end{center}
	\vspace{-0.43cm}
\end{table}

In Figure \ref{fig:simulations}, we show the estimation of $a_1$ and $b_1$ in a single run. We can appreciate that when $r=1000$, $a_1$ and $b_1$ are generally well approximated  and the abrupt changes are promptly detected. Some picks are visible (for instance, immediately after 0.4 ms in the $a_1$ graph), due to locally large noise (some oscillations were highlighted also in \cite{li11}). In the case $r=100$, instead the changes are less promptly detected (in particular, for $a_1$). However, a good approximation is achieved after some steps. Moreover, a smaller number of iterations makes the algorithm more conservative, thus less affected by locally large noise: for example, the pick after 0.4 ms is canceled.

In Figure \ref{fig:support}, we show the signed error $\widehat{v}_{i,sm}-v_{i,sm}$, $s=1,\dots,66$. The values 0 and 10 on the $x$-axis respectively represent $a_t$ and $b_t$; the others are the null ones.  Green circles and blue squares represent the estimations, at each $sm$, respectively for $r=100$ and $r=1000$. For $r=1000$,  $\widehat{v}_{i,sm}-v_{i,sm}$ is smaller on the null parameters if compared to non-null parameters $a_1$ and $b_1$. This means that the support (\emph{i.e.}, the positions of the non-null components in the vector $v_t$) is well detected.
\section{CONCLUSIONS}
In this work, we have presented a method for online optimization for time-varying sparse problems. We have considered an Elastic-net cost functional and proposed an online IST algorithm. We have then proved that this method is successful in terms of dynamic regret. Moreover, we have shown how to apply the algorithm to a problem of recursive identification for a TVARX model, and we have presented a few numerical results. In future work, other algorithms (\emph{e.g.}, ADMM) will be considered for this problem and  other sparse  (possibly not strongly convex) models will be investigated.

%%%%%%%%%%%%%%%%%%%%%%%%%%%%%%%%%%%%%%%%%%%%%%%%%%%%%%%%%%%%%%%%%%%%%%%%%%%%%%%%
%\section{ACKNOWLEDGMENTS}
%The author wishes to thank Diego Regruto and Vito Cerone for useful discussions and suggestions.

%%%%%%%%%%%%%%%%%%%%%%%%%%%%%%%%%%%%%%%%%%%%%%%%%%%%%%%%%%%%%%%%%%%%%%%%%%%%%%%%

\bibliographystyle{plain}
\bibliography{refs}

\end{document}